\documentclass{article}
\usepackage{geometry}
\usepackage{chips}
\usepackage{dyck}
\usepackage{paper}

\begin{document}
\begin{center}
  {\Large \textbf{Enumerating Pyranim Isomorphisms}}\\[0.6em]
  Yoav Binyamin$^{1}$, Evan Kniffen$^{1}$\\[0.25em]
  {\small
  $^{1}$Department of Mathematics, Texas A\&M University
  }\\[0.8em]
\end{center}
\begin{abstract}
We study a Nim-adjacent poset game, where chips are arranged in a pyramid and two players take turns removing chips that do not currently support chips to the upper-left and upper-right; upon removing this chip, all supported chips above are also removed. Play continues until the final chip is removed. We use a bijection between the game states, Dyck paths, and unit interval graphs to count the number of possible game states up to isomorphism.
\end{abstract}

\begin{center}
\end{center}

\vspace{.25cm}

\tableofcontents

\newpage
\section{Introduction}
\textbf{Pyramid Nim} is defined in \cite{pyranim}, studied as a combinatorial game played on a directed acyclic graph. However, the authors (i.e. Binyamin) came up with \textbf{Pyranim} independently in its more specific play on triangular starting positions. Historically, both formulations independently used ``pyramid" rather than ``triangle" to convey the stability of the positions in the game.

Our formulation of this 2-player impartial Nim game Pyranim is as follows: you start with an $n$-tier triangle of chips (consisting of $\tfrac{n(n+1)}2$ chips total). Player 1 goes first, and can take any of the chips on the outer layer of the pyramid. By ``outer layer," we mean every chip that does not support two chips on the row directly above it; that is, every chip which can be seen observing the game from above. Upon doing so, it will also remove all chips supported by the removed chip. Player 2 does the same with the modified pyramid and play continues until the final chip is taken, and the last to take a chip is the winner. We will continue with an example game to describe the proceedings:

\begin{example} For a 3-tier game with six total chips, the game starts with Player 1 being able to remove any chip other than the middle chip on the third layer. 
    \[\text{\PY{UUUUDDDD}}\]
\begin{enumerate}
    \item At this step, Player 1 \textit{cannot} remove the bottom middle chip, but \textit{can} remove any other chip. Player 1 removes the left chip on the second tier, taking the apex chip with it. This leaves the three bottom chips and the right chip on the second tier.
    \[\hspace{-1cm}\text{\PY{UUDUUDDD}}\]
    \item Player 2 removes the lower right chip, leaving two chips on the lowest row.
    \[\hspace{-1cm}\text{\PY{UUDUDDUD}}\]
    \item Player 1 removes one chip, leaving the final chip to Player 2, who wins the game.
    \[\hspace{-1cm}\text{\PY{UUDDUDUD}}\]
\end{enumerate}
\end{example}

\subsection{Setup}
To formally define the game, we will treat each chip as a discrete interval containing the chips supporting it. Let \(\mathcal{P}_n := \{[i, j] : 1 \leq i \leq j \leq n\}\) be ordered by containment: \[[a, b] \leq [i, j] \iff i \leq a \leq b \leq j.\] Write $c_{i,j} = [i, j]$. The bottom chips are $c_{i,i}$; the immediate
supports of $c_{i,j}, \, i < j$, are $c_{i,j-1}$ and $c_{i+1 ,j}$. This structure will allow us to wield a multitude of wonderful algebraic and combinatorial tools.

\begin{definition}
    A \textbf{stable position} $X$ is an order ideal in $\mathcal{P}_n$ under containment:
\[c_{i,j} \in X,\quad i \leq a \leq b \leq j \implies c_{a,b} \in X.\]
Let \(\mathcal{T}_n\) be all stable positions and \(\mathcal{S}_n\) those with \(c_{1,1},\dots,c_{n,n}\in X\). Set \(\mathcal{T}_0=\mathcal{S}_0=\{\emptyset\}\).
\end{definition}
 We interpret $\mathcal{T}_n$ as the set of all of the legal games generated by normal play and \(\mathcal{S}_n\) as the set of all of the legal games with a full bottom row of chips (\(\mathcal{S}_n\subseteq \mathcal{T}_n\)). This is an important distinction when considering raw states generated by legal play on different pyramid sizes, since \[\PY{UUDD}\hspace{2em}\raisebox{2em}{$\not=$}\hspace{2em}\PY{UDUUDD}.\]

\begin{definition}
    A surviving chip $c_{i,j}\in X$ is \textbf{visible} if its two possible immediate upper chips $c_{i-1,j}$
and $c_{i,j+1}$ are not both present (indices outside $\mathcal{P}_n$ are absent).

Only visible chips are playable.
Removing one chip forces the recursive removal of every chip that subsequently lacks a support. That is, if $s(x)$ is the maximal order ideal contained in $x \subseteq \mathcal{P}_n$, then the state that forms from removing $c_{i, j}$ is $s(X \setminus c_{i, j} )$.
Let the option set, $\opt(X)$, be the set of distinct successor positions. That is,
\begin{equation*}
    \opt(X) = \{ s(X \setminus c_{i, j}) : c_{i, j} \text{ is visible}\}.
\end{equation*}
\end{definition}

\subsection{Isomorphisms}
While some states may be \textit{literally} different as defined, they play out the same \textit{in practice} (an exmaple of which can be seen in the above figure). For example, reflections and translations of a game play out the exact same. Inspired to reduce the number of effective game states, we recursively define isomorphism.
\begin{definition}\label{recursive}
    Let $X$ and $Y$ be two states. Then $X$ and $Y$ are \textbf{isomorphic} ($X \simeq Y$) as games, if, given that the next legal moves from $X$ are
    $\opt(X)=\{X_1, \dots, X_m\}$ and similarly for $Y$, then there is a bijection (called \textbf{isomorphism})
    \[\varphi: \opt(X) \longrightarrow \opt(Y)\quad\text{ such that }\quad\forall i \in [m], \,\varphi(X_i) \simeq X_i.\] Note that the empty state of 0 chips  is isomorphic only to itself. Game isomorphism will sometimes be denoted \(\simeq_{game}\) to differentiate with graph isomorphisms when unclear.
\end{definition}

    \begin{example}
        \[\PY{UUDDUUUDDD}\hspace{1em}\raisebox{5em}{$\simeq$}\hspace{1em}\PY{UDUUDUUDDD}\]\hspace{.5em}\[\PY{UUDUUDDD}\hspace{1em}\raisebox{4em}{$\simeq$}\hspace{1em}\PY{UUUDDUDD}\hspace{5em}\PY{UUDUDUDD}\hspace{1em}\raisebox{4em}{$\not\simeq$}\hspace{1em}\PY{UUUDDDUD}
    \]
    \end{example}
As expected, this isomorphism is an equivalence relation. Therefore, for all \(n\geq0\), put \[\mathcal{K}_n=\mathcal{T}_n/\simeq_{game}\qquad \text{and}\qquad K_n=|\mathcal{K}_n|.\] \(K_n\) represents the number of states on each \(n\)-pyramid up to isomorphism (with $K_0=1$). The rest of the paper focuses on finding a closed form of $K_n$ and its ordinary generating function \begin{equation}\label{gf}
    K(x)=\sum_{n\geq0}K_nx^n.
\end{equation}

\section{Structural Bijections}
    In this section we will discuss the bijections present under normal state-counting as well as modulo isomorphism. In particular, the former bijects to Dyck paths and the latter bijects into unlabeled unit interval graphs.
\subsection {Dyck Paths \& Words}
    Recall that a Dyck path is a north-east lattice path from $(0,0)$ to $(n,n)$ that never goes strictly below the line $y=x$ (i.e. subdiagonal). The set of Dyck paths of length $2n$, and thus semilength $n$, will be denoted $\dyck_{2n}$. Dyck paths can also be written as words of the Dyck language, with up-steps written as $U$ and down-steps written as $D$. We force that every prefix has at least as many $U$'s as $D$'s for a valid Dyck word. \hspace{.5em}
    
    \[\DyckPathWord{UUDUDD}\]\begin{center}
        $UUDUDD$
    \end{center} 
    \begin{definition}
        A Dyck path is \textbf{primitive} if and only if it returns to the diagonal at only the final down-step. This can be extended to words---a word is primitive if no proper nonempty prefix is itself balanced in the number of $U$'s and $D$'s. Alternatively, removing the first $U$ and last $D$ results in another Dyck word.
    \end{definition}  
    We will use a well-known combinatorial result,\begin{equation}\label{catalan}
        |\dyck_{2n}|=C_n=\frac{1}{n+1}\binom{2n}{n},
    \end{equation}
    where $C_n$ denotes the Catalan numbers \cite{stanley-catalan}.
    \subsection{Dyck Correspondence}
    Intuitively, we can see that our pyramids fit nicely into the full boxes generated by the area between the Dyck path and the diagonal. 
    \begin{prop}
    Each stable position (i.e. legal state) for an $n$-pyramid game has a unique Dyck path of semilength $n+1$ and corresponding Dyck word of length $2(n+1)$. This process is reversible, inducing a clear bijection. Therefore by Equation \ref{catalan}, \[|\mathcal{T}_n|=C_{n+1}=\frac{1}{n+2}\binom{2n+2}{n+1}.\]
\end{prop}

\hspace{.5em}
    \begin{center}
\PY{UUDDUUUDDD}
\hspace{2em}\raisebox{5em}{$\longleftrightarrow$}\hspace{3em}
\DyckPathWord{UUDDUUUDDD}
\end{center}

\medskip
Both of these are represented by the Dyck word $UUDDUUUDDD$. This longer word can also be written in the form $U^2D^2U^3D^3$. A full [resp. empty] \(n\)-pyramid can be written as $U^{n+1}D^{n+1}$ [resp. $(UD)^{n+1}$].
\subsection{Unit Interval Graphs}

We now show that \(\mathcal{K}_n\) corresponds to a subclass of known graphs. 

 \begin{definition}
     A \textbf{unit interval graph} is a graph formed from a set of intervals on the real line, where each interval has a length of exactly $1$. Each vertex represents an interval, and an edge connects two vertices if and only if their corresponding intervals overlap.
 \end{definition}
 \begin{example}

 \end{example}
\hspace{1.2cm}\begin{tikzpicture}[scale=1.5, every node/.style={circle, draw, minimum size=0.6cm, inner sep=0pt}]
        
        % Drawing Intervals
        \begin{scope}[shift={(0,1.5)}, every node/.style={}]
            \draw[gray, fill=gray!20] (0, 0.6) rectangle (1, 0.7) node[right, black] {$A: [0, 1]$};
            \draw[gray, fill=gray!20] (0.4, 0.3) rectangle (1.4, 0.4) node[right, black] {$B: [0.4, 1.4]$};
            \draw[gray, fill=gray!20] (0.9, 0.0) rectangle (1.9, 0.1) node[right, black] {$C: [0.9, 1.9]$};
            \draw[gray, fill=gray!20] (1.5, -0.3) rectangle (2.5, -0.2) node[right, black] {$D: [1.5, 2.5]$};
            % Axis
            \draw[->] (-0.2, -0.6) -- (2.8, -0.6) node[right] {$\mathbb{R}$};
            \foreach \x in {0,1,2} \draw (\x, -0.55) -- (\x, -0.65) node[below] {\x};
        \end{scope}
        
        % Drawing Labeled Graph Vertices
        \node (A) at (5,1.2) {A};
        \node (B) at (6,1.2) {B};
        \node (C) at (7,1.2) {C};
        \node (D) at (8,1.2) {D};
        
        % Drawing Edges based on overlaps
        \path[-] (A) edge (B);
        \path[-] (B) edge (C);
        \path[-] (C) edge (D);
        \path[-] (A) edge[color=black, bend left=40] (C); % A and C intersect between 0.9 and 1.0
        
    \end{tikzpicture}

 The following characterizations are equivalent to unit interval graphs: indifference graphs, proper interval graphs, and $K_{1,3}$-free interval graphs \cite{BogartWest1999, IntervalRoberts}. We will use these terms interchangeably depending on what property we are trying to highlight. If the labels $A$, $B$, $C$, and $D$ are removed, we get an \textbf{unlabeled unit interval graph} (or UUIG, for brevity). Permuting the labels of a labeled unit interval graph creates an isomorphic labeled unit interval graph, but these two graphs are \textit{equal} as UUIGs.

\subsection{Unlabeled Unit Interval Graph Correspondence}
\begin{definition}
    Construct the graph $\Gamma(X)=(V,E)$ by \[V=\{i:c_{i,i}\in X\}\qquad \text{and}\qquad E=\{\{i,j\}:i<j,\, c_{i,j}\in X \}.\] 
\end{definition}
Therefore,
\begin{equation}
    |X|=|V(\Gamma(X))|+|E(\Gamma(X))|.
\end{equation}

\begin{remark}
    We interpret this map as follows: map each game state $X$ to its Dyck word. Remove the first U and last D.  Let interval $i$ be the interval with the $i$th $U$ as the left endpoint and the $i$th $D$ as the right endpoint, interpreting intervals where the right endpoint appears before the left endpoint as missing. Then $\Gamma (X)$ is the proper interval graph of these intervals.
\end{remark}

\begin{lemma}\label{umbrella}
    $\Gamma(X)$ has the umbrella property
\[i < j < k,\, \{i,k\} \in E \implies \{i,j\}, \{j,k\} \in E.\]
$\Gamma(X)$ is therefore a proper interval graph; equivalently, a unit interval graph. Every unit interval
graph with a chosen proper interval ordering $1<\cdots<n$ yields one supported ideal \cite{Kato2024}: retain each $c_{i,i}$ and retain $c_{i,j}$ if and only if $\{i,j\}\in E$. 
\end{lemma}
\begin{proof} If $c_{i,k}$ survives, all its subintervals, including $c_{i,j}$ and $c_{j,k}$, survive. Conversely, the umbrella
property, iterated over indices within $[i, j]$, makes every pair inside an edge $\{i,j\}$ an edge as well. The
chips specified by an ordered graph are therefore an ideal. The proper–unit equivalence is proved in \cite{BogartWest1999}.
\end{proof}

We can pass from \textit{labeled} unit interval graphs to \textit{unlabeled} unit interval graphs by looking at these constructions modulo isomorphism.
\begin{thm}[Isomorphism Equivalence]\label{targetthm}
    For \(X,Y\in\mathcal{T}_n,\) \[ \Gamma(X)\simeq\Gamma(Y) \iff X\simeq_{game}Y.\]
\end{thm}
We will spend the majority of the paper proving this theorem. The benefits of this will become clear after taking note of known results regarding UUIGs. 
\begin{cor}\label{form}
Let \[\mathcal{S}^*_n=\mathcal{S}_n/\simeq_{game}\]

    By Theorem \ref{targetthm}, $\mathcal{S}_n^*$ bijects with UUIGs on $n$ vertices (thus the two objects are equinumerous). Let $u_n$ denote the counts of those graphs---corresponding to the sequence A005217 \cite{OEISA005217}. Then $|\mathcal{S}_n^*|=u_n$.
As consequence of working modulo game isomorphism, we get 
%\begin{equation}
%    \mathcal{K}_n=\bigcup_{k=0}^n\mathcal{S}^*_k.
%\end{equation}
%Taking the cardinality of both sides to get the %desired count,
\begin{equation}\label{target}
    K_n=\sum_{k=0}^nu_k.
\end{equation}

Indeed, we can consider states of different-size support separately; if a state $X$ has $k$ bottom chips, then one can trivially find an isomorphism between $X$ and an element of $\mathcal S_k$.
\end{cor}

\section{Isomorphism Equivalence}
\subsection{The Forward Implication}
\begin{lemma}\label{cascade}
    For all visible $c_{i,j}$ in the order ideal $X$, removing that chip results in the position \[M_{c_{i,j}}(X)=X\setminus\{c_{a,b}\in X : a\leq i \leq j \leq b\}.\] Where distinct playable chips give distinct legal successor positions.
\end{lemma}
\begin{proof}
    Every superinterval of $[i, j]$ requires it as support through a chain of immediate supports and
must disappear. Deleting the entire set of superintervals in $M_{c_{i,j}}(X)$ leaves an ideal: any subinterval of a surviving chip still survives, as it is unrelated to the original chip---or any chip removed thereafter. Thus no other chip can lose a support. To see that distinct chips $c_{i,j}$ and $c_{k,l}$ yield distinct positions, assume without loss of generality that $[k,l]\not\subseteq [i,j]$. Therefore, the interval $[i,j]$ is \textit{not} a superinterval of $[k,l]$, meaning $c_{i,j}$ is preserved when $c_{k,l}$ is removed. Because $c_{i,j}$ is eliminated when being played, we have \[c_{i,j}\in M_{c_{k,l}}(X)\setminus M_{c_{i,j}}(X),\] so the successor positions are distinct.
\end{proof}
\begin{lemma}\label{inter}
    Graph components occur in consecutive runs of bottom
indices and no retained chip spans two components. A chip move modifies only its own component.
\end{lemma}
\begin{proof}
    If consecutive vertices lack an edge, a crossing edge would force that edge by Lemma \ref{umbrella}, so there are no crossing edges. If a bottom chip is missing, no chip containing that
singleton can survive. A retained interval $[i, j]$ forces all consecutive edges between $i$ and $j$, and hence
cannot span distinct components. Therefore, \(M_{c_{i,j}}(X)\) is restricted to one component.
\end{proof} 

\begin{lemma}\label{reversal}
    Let $G, H$ be two connected proper unit interval graphs, represented by intervals $(I_i)_{i=1}^m$ and $(J_i)_{i=1}^m$, respectively, where the intervals are listed in order based on their left endpoint from left to right. Then $G \simeq H$ if and only if
    $$\forall i, j \in [m], I_i \cap I_j = \emptyset \iff J_i \cap J_j = \emptyset,$$

    or
    $$\forall i, j \in [m], I_i \cap I_j = \emptyset \iff J_{m+1-i} \cap J_{m+1-j} = \emptyset.$$
\end{lemma}
This result is the same as Corollary 2.5 from \cite{Deng}.

\begin{prop}[The Forward Implication]\label{forward}
 For \(X,Y\in\mathcal{T}_n,\)\[\Gamma(X)\simeq\Gamma(Y)\implies X\simeq_{game}Y.\]    
\end{prop}
\begin{proof}

By Lemma \ref{inter}, \(\Gamma(X)\)'s connected components' vertices occupy consecutive indices, and no
surviving chip crosses from one connected component to another. Since
$\Gamma(X)\simeq\Gamma(Y)$, we may match the connected components of $X$ with those of
$Y$. Lemma \ref{reversal} says that the ordered chip patterns in a matched pair agree, possibly after reversing one
connected component. Thus we obtain a bijection $F:X\to Y$ by moving entire connected components
and, if necessary, reversing their local indices. On a connected component of
length $m$, reversal is defined by the map
\[
    c_{ij}\longmapsto c_{m+1-j,\,m+1-i},
\]
where $i,j$ are interpreted as local indices.

We claim that this same map works throughout the game, even after
bottom chips have been removed. For any stable subposition $Z\subseteq X$,
put $\widehat F(Z)=\{F(c):c\in Z\}$. Shifting or reversing a connected component
preserves containment of intervals (since chips in different connected components cannot contain one another). Hence $\widehat F(Z)$ is stable. Moreover, the two
possible chips immediately above $c$ are sent to the two possible chips
immediately above $F(c)$, with their roles exchanged under reversal.
Note that any upper chip crossing a connected component boundary is absent on both sides.
Consequently, \[c \text{ is visible in } Z \iff F(c) \text{ is visible in } \widehat F(Z).\]

The cascade after playing $c$ removes precisely the surviving chips
whose intervals contain $c$. Since $F$ preserves containment,
\[
    \widehat F\bigl(M_c(Z)\bigr)
       =M_{F(c)}\bigl(\widehat F(Z)\bigr).
\]
Thus by Lemma \ref{cascade}, $\widehat F$ matches every legal option of $Z$ with a legal option
of $\widehat F(Z)$, and the inverse connected component map matches them back.
Because $\widehat F$ is injective on positions, distinct actual
successors remain distinct, even when they belong to the same
game-isomorphism class. We keep this one connected component map for all descendants;
we do not compress and relabel each successor anew.

Starting with $\widehat F(X)=Y$, the displayed identity gives a
bijection of options at every subsequent position. Induction on the
number of surviving chips therefore gives $X\simeq_{\mathrm{game}}Y$.
\end{proof}
\begin{remark}
    The above proposition relies on the known result  that graphs are isomorphic if and only if they have isomorphic connected components, as well as the result that connected unit interval graphs are isomorphic if and only if they were created from the same or reversed order of intervals (Lemma \ref{reversal}). In the game sense, this corresponds to translations of connected components, permutations of the order of appearance of connected components, and reversals within connected components. However, these clearly preserve isomorphism. The converse implication, however, says that these are in fact the \textit{only} ways for two game states to be isomorphic; we prove the converse in Section \ref{section6}.
\end{remark}

\subsection{Reconstruction}
The converse implication of Theorem \ref{targetthm} means that for a state $X$, knowing the graphs of $\text{Opt}(X)$, we can reconstruct the original $\Gamma (X)$ uniquely. However, we will prove an even \textit{stronger} statement: rather than the entirety of $\text{Opt}(X)$, one only needs the states where a \textit{bottom chip} was removed from $X$. This is captured by the following definition and theorem.

\begin{definition}
    For a unit interval graph $G$, choose a proper interval ordering within each connected component. Select its first and last vertices if it has at least two vertices, or its unique vertex if it is a singleton. The \textbf{component-end deck} is the multiset
    $$\partial G = \{ \{ [G - v] : v \text{ is a selected component end} \} \},$$

    where $[ \cdot ]$ denotes graph isomorphism class and $G - v$ is the graph gotten by removing $v$ and all of the edges connected to it. An element of $\partial G$ is called a \textbf{card}.
\end{definition}

\begin{thm}[Reconstruction]\label{reconstruction}
    Every unit interval graph $G$ is determined up to isomorphism by
    \begin{equation}\label{triple}
        (|V (G)|, |E(G)|, \,\partial G).
    \end{equation}
\end{thm}

\subsubsection{Connected Reconstruction}
We'll first prove the theorem for connected $G = \Gamma (X)$, and then use the result for each of the connected components of a disconnected $G$. When $G$ is connected, $X$'s Dyck word is primitive.

\begin{definition}\label{lrpd}
    Let $w$ be a primitive Dyck word, and write $w = U^r D z$ for $z$ a word on $U$'s and $D$'s. Then
    \begin{enumerate}[label=(\roman*)]
        \item $\rho(w)$ flips the order of the $U$ and $D$, and then swaps every $U$ to a $D$, and vice versa.
        \item $[w]  := \{ w, \rho(w) \}$.
        \item $L(w) := U^{r-1} z$ deletes the first $U$ and first $D$.
        \item $R(w) := \rho (L (\rho (w)))$ deletes the last $U$ and last $D$.
        \item For $0 \leq d \leq r$,
        $P_d (w) := U^{d+1} D U ^{r-d} D z.$
    \end{enumerate}
\end{definition}
\begin{remark}
    We interpret this as follows:
    \begin{enumerate}[label=(\roman*)]
        \item $\rho$ flips the game. If one player sees $w$, then the player on the other end of the table sees $\rho(w)$.
        \item If the proper interval graph of $w$ is connected, then $[w]$ is the set of all primitive Dyck words with the same interval graph as $w$. Indeed, by Lemma \ref{reversal}, the word and its reversal are the only possible representatives.
        \item $L$ removes the leftmost (bottom) chip.
        \item $R$ removes the rightmost (bottom) chip.
        \item $P_d$ adds $d$ chips on the left. The original $w$ had $r-1$ chips to add on top of, so $d \leq (r-1)+1 = r.$
    \end{enumerate}
\end{remark}

\begin{thm}[Conjugacy]\label{conjugation}
    \cite{KarhumakiNotes} Let $A, B, z$ be words and $|A| = |B| = h > 0$. Then \[Az = zB\] if and only if for some $t \geq 0$ and words $p, q$ with $|p| < h$,
    \begin{equation*}
        A = pq, \hspace{0.5cm} B = qp, \hspace{0.5cm} z = A^t p.
    \end{equation*}
\end{thm}
\begin{proof}
    If $|z| \geq h$, then comparing the initial $h$ letters, $z = Az_1$ for some word $z_1$, where $|z_1| = |z| - h$. Canceling, $Az_1 = z_1 B$. Continuing inductively, the process terminates because $z$ has finite length. We get $z = A^t p$ and $Ap = pB$ for some word $p$, with $|p| < h$. Inspecting the first $h$ letters, $A = pq$ for some word $q$, and $B = qp$. The converse follows computationally from $(pq)^{t+1} p = (pq)^t p (q p)$.
\end{proof}

\begin{lemma}[Nine-case Orientation]\label{nine-case}
    Let $w$ be primitive, $d \geq 1$, and assume both $P_d(w)$ and $P_d(\rho (w))$ are defined. Then
    $$ [ R(P_d(w)) ] = [R (P_d( \rho (w) ) ] \implies w = \rho(w).$$

    In particular,
    $$ [ R(P_d(w)) ] = [R (P_d( \rho (w) ) ] \implies [P_d (w)] = [P_d (\rho (w))].$$
\end{lemma}
\begin{proof}
    If $w = U^k D^k$, then it is self-reflecting, so the result is trivially true. Otherwise, write for some word $z$
    $$w = U^r D z U D^s, \hspace{0.5cm} \# U(z) - \# D(z) = s-r.$$

    Then
    \begin{equation}\label{rho_w}
        \rho(w) = U^s D \rho(z) U D^r.
    \end{equation}

    Put $\alpha = r-d \geq 0$ and $\beta = s - d \geq 0$. Then
    $$a = R( P_d (w)) = U^{d+1} D U^{\alpha} D z D^{s-1},$$
    $$b = R(P_d(\rho(w))) = U^{d+1} D U^\beta D \rho(z) D^{r-1}.$$

    By assumption, $a = b$ or $a = \rho(b)$. If $a = b$ then $\alpha = \beta$, so $r=s$ and $z = \rho(z)$. Therefore, $w = \rho(w)$, and we are done.

    If $a = \rho(b)$, then
    \begin{equation}\label{word_ninecase}
        U^{d+1} D U^\alpha D z D^{s-1} = U^{r-1} z U D^\beta U D^{d+1}.
    \end{equation}

    %Also, since $b = \rho(a)$, then we have
    %$$U^{d+1} D U^{\beta} D z D^{s-1} = U^{r-1} z U D^\alpha U D^{d+1}. $$

    Inspecting Equation \ref{word_ninecase}, $d+1 \geq r-1$ and $d+1 \geq s-1$. Hence, $0 \leq \alpha \leq 2, 0 \leq \beta \leq 2.$ This gives us the titular nine cases.

    Let $A_\alpha = U^{2 - \alpha} D U^\alpha D$ and $B_\beta = U D^\beta U D^{2- \beta}.$ Then $A_\alpha z = z B_\beta$. Using Theorem \ref{conjugation}, we arrive at the following conclusions as to what $z$ is for each pair $(\alpha, \beta)$, where $t \geq 0$:

    \begin{center}
    \begin{tabular}{|>{\centering\arraybackslash}m{0.08\textwidth}|
  *{3}{>{\centering\arraybackslash}m{0.16\textwidth}|}}\hline
    & $\beta = 0$ & $\beta = 1$ & $\beta = 2$ \\\hline
    $\alpha = 0$ & $***$ & $*$ & $**$ \\
    & $(UUDD)^t$ & N/A & $(UUDD)^t U$\\\hline
    $\alpha =1$ & $*$ & $***$ & $*$\\
    & N/A & $(UD)^t$ & N/A \\\hline
    $\alpha = 2$ & $**$ & $*$ & $***$\\
    & $(DUUD)^t D$ & N/A & $(DUUD)^t DU $ \\\hline
    
    \end{tabular}
    \end{center}

    \begin{itemize}
        \item For the cases marked with (*), one \textit{cannot} find $p, q$ for which $|p| < 4$, $A_\alpha = pq,$ and $B_\beta = qp$. This is checked by noting that the $4$ possible $p$ do not satisfy the relations. Therefore, these cases are impossible.
        \item For the cases marked with (**), one \textit{can} find $p, q$ for which $|p| < 4$, $A_\alpha = pq,$ and $B_\beta = qp$. However, the number of $U$'s and $D$'s for $w$ is not the same, so these cases are also impossible. Indeed,
            \begin{enumerate}[label=(\roman*)]
                \item $(0, 2)$: $p = U, q = UDD$ is the only option that works. \[\# U(z) - \#D(z) = s-r = \beta - \alpha = 2.\] However, $z$ clearly has $\# U(z) - \#D(z) = 1 \neq 2$.
                \item $(2, 0)$: $p = D, q = UUD$ is the only option that works. \[\# U(z) - \#D(z) = s-r = \beta - \alpha = -2.\] However, $z$ clearly has $\# U(z) - \#D(z) = -1 \neq -2$.
            \end{enumerate}
        \item The cases marked with (***) \textit{are} possible. One can find $p, q$ for which $|p| < 4$, $A_\alpha = pq,$ and $B_\beta = qp$. Then, these choices imply that $r=s$, since $\alpha = \beta$. They also imply that $z = \rho(z)$. Indeed,
            \begin{enumerate}[label=(\roman*)]
                \item $(0, 0)$: We can choose any $p$ of the $4$, but only $p = \emptyset,\, q = UUDD$ satisfies the condtion \[\# U(z) - \#D(z) = \beta - \alpha = 0.\] Therefore, $z = (UUDD)^t \emptyset = (UUDD)^t$, so $z = \rho(z)$.
                \item $(1, 1)$: We can choose either $p= UD,\, q = UD$ or $p = \emptyset,\, q = UDUD$. In the former, \[z = (UDUD)^t UD,\]so there are \textit{odd} many $UD$. Meanwhile, in the latter, \[z = (UDUD)^t \emptyset = (UDUD)^t,\]so there are \textit{even} many $UD$. Combining the two cases, $z = (UD)^t$ for some $t \geq 0$. This satisfies \[\# U(z) - \#D(z) = \beta - \alpha = 0.\] Finally, $z = \rho(z)$.
                \item $(2, 2)$: $p = DU, q = UD$ is the only choice. $z = (DUUD)^t DU$, \[\# U(z) - \#D(z) = \beta - \alpha = 0,\] and $z = \rho(z)$.
            \end{enumerate}

    \end{itemize}

    Therefore, the only possible cases satisfy $r=s$ and $\rho(z) = z$. By inspecting Equation \ref{rho_w}, $w = \rho(w)$. The ``in particular" follows immediately, since $P_d(w) = P_d(\rho(w))$ from Definition \ref{lrpd}.

\end{proof}
\begin{prop}[Connected Reconstruction]\label{connected_reconstruction}
    Every connected unit interval graph $G$ is determined up to isomorphism by
    $$(|V(G)|, |E(G)|, \,\partial G).$$

\end{prop}
\begin{proof}
The cases of $|V(G)| \leq 2$ are trivially true. Suppose $|V(G)| \geq 3$, and $G,H$ are connected unit interval graphs with equal triples. Select a card $[A]$ from $\partial G$. Reflect either parent graph, if needed, so that $[A]$ is its left-end card. Since consecutive vertices in a connected proper interval ordering are adjacent, deleting an end preserves connectedness. Give $A$ a primitive word $w$ and let $d$ be the number of chips removed. The degree of the deleted vertex is the same in both parents:
\[
d  - 1= |E(G)| - |E(A) |= |E(H)| - |E(A)|.
\]
Of course, $d \geq 1$. The deleted left-end vertex contributes the first $U$ and first $D$ of the parent word, with that $D$ following exactly $d$ of the remaining $U$'s. Deleting these two letters recovers either $w$ or $\rho(w)$.

Therefore, each parent has a word which is either $P_d(w)$ or $P_d(\rho(w))$. If the operation is defined for only one of $w$ or $\rho(w)$, then the two parents $G,H$ are represented by the same word, and are therefore in the same isomorphism class. If the operation is defined for both $w$ and $\rho(w)$, then they are either both represented by the same word, or one of $G,H$ is represented by $P_d(w)$ and the other by $P_d(\rho(w))$.

In the latter case, remove one occurrence of $[A]$ from each of the equal two-entry end decks. The remaining entries are the right-end cards, so
\[
[R(P_d(w))] = [R(P_d(\rho(w)))].
\]
This cancellation is valid even when both entries of an end deck equal $[A]$. Hence, by Lemma~3.8,
\[
[P_d(w)] = [P_d(\rho(w))],
\]
so that $G,H$ are isomorphic.
\end{proof}

\subsubsection{Disconnected Reconstruction}
We'll prove Theorem \ref{reconstruction}.

\begin{proof}
We argue by strong induction on $|V(G)|$. The cases of $|V(G)| \leq 2$ are determined by the vertex and edge counts. If $G$ has $q$ nonsingleton components and $s$ singleton components, then
\[
|\partial G| = 2q+s.
\]
Thus, for $|V(G)| \geq 3$, a connected graph has two cards, whereas a disconnected graph has at least three. The triple therefore determines whether $G$ is connected.

Since the statement is proven for connected $G$ by Proposition \ref{connected_reconstruction}, let $G$ be disconnected. If $G$ is edgeless, this is trivially true. Otherwise, let $M$ be the order of the largest connected component of $G$, in which case $M \geq 2$. Removing an end of a connected proper interval graph preserves connectedness and reduces its order by $1$ if it is not a singleton, and deletes it if it is a singleton. Indeed, consecutive vertices in a connected proper interval ordering are adjacent; otherwise, Lemma \ref{umbrella} would forbid every edge crossing their gap.

Since $G$ has at least two connected components, there is a card of $G$ for which a component of order $M$ is maintained. Consequently, $M$ is the largest component order occurring in any card of $\partial G$.

For the following, let $T$ denote an isomorphism class of connected proper interval graphs on $M$ vertices. For a card $F \in \partial G$, let $b_T(F)$ be the multiplicity of $T$ as a component of $F$. Write the vector
\[
b(F) = (b_T(F))_T,
\]
indexed by $T$.

If $b(F) = b(F')$ for every $F,F' \in \partial G$, then every component is in the same class $T$. Indeed, a smaller component would allow a deletion preserving every order-$M$ component, while deletion in an order-$M$ component decreases one coordinate. Two different order-$M$ types would likewise allow deletions decreasing different coordinates. Either possibility would give unequal vectors. Since $G$ has at least two components, a copy of $T$ survives in every card, so its type is determined. Thus, $[G]$ is the class of $|V(G)|/M$ copies of $T$.

Otherwise, if the vectors $b(F)$ are not all equal, then for each $T$ there is an original component other than a copy of $T$. Deleting an end in that component retains all copies of $T$, and no deletion creates a new component of order $M$. Therefore, for each $T$, the number of components of class $T$ in $G$ is
\[
m_T := \max_{F \in \partial G} b_T(F).
\]

Therefore, we know all of the components of $G$ of order $M$ up to isomorphism. Now let $Q$ be their disjoint union, and consider the subgraph $R = G-Q$ made by excluding these components. If $R$ is empty, then $G=Q$ is already determined.

Otherwise, we know $\partial R$ exactly from the cards $F \in \partial G$ with
\[
b_T(F) = m_T \qquad  \forall T.
\]
Indeed, deletion in an order-$M$ component decreases one of these coordinates, while deletion in a smaller component leaves all of them unchanged. Thus, these cards are precisely the end vertex removals in components with fewer than $M$ vertices. Removing $Q$ from each selected card, while preserving every occurrence, gives the multiset $\partial R$. Here, $F-Q$ means the class $F$ with all of its components of order $M$ removed.

We also know
\[
|V(R)| = |V(G)| - |V(Q)|,
\qquad
|E(R)| = |E(G)| - |E(Q)|.
\]
Hence, the entire triple (\ref{triple}) for $R$ is determined. Since $R$ has fewer vertices than $G$, the inductive hypothesis determines $R$ up to isomorphism. Knowing all of the other components, $G$ is determined up to isomorphism.
\end{proof}

\subsection{The Converse Implication}\label{section6}
%\begin{lemma}[Invariants]\label{invariants}
    We will now define some invariants that occur between games. Let $X \simeq_{\text{game}} Y$. Then the number of chips at the bottom (labelled $c_{ii}$) and the total number of chips of $X$ and $Y$ are the same. Additionally, if $\ell_{\max} (X)$ is the maximum number of moves possible on $X$ before the null position, and $\ell_{\min} (X)$ is the minimum number of moves possible on $X$ before the null position, then
    \begin{enumerate}[label=(\roman*)]
        \item $\ell_{\max} (X) = |V (\Gamma(X))| + |E(\Gamma (X)| = |X|$ is the total number of chips of $X$, and
        \item $\ell_{\min} (X) = |V(\Gamma(X))|$ is the total number of bottom chips of $X$.
    \end{enumerate}

    Consequently, $$\ell_{\max} (X) = \ell_{\max} (Y) , \qquad\qquad \ell_{\min} (X) = \ell_{\min} (Y),$$  $$|V(\Gamma(X))| = |V(\Gamma(Y))|, \quad\quad |E(\Gamma(X))| = |E(\Gamma(Y))|.$$ 
%\end{lemma}
%\begin{proof}
    (Indeed, the invariance of $\ell_{\max}$ is a trivial inductive argument on the number of chips.
    
    The invariance of $\ell_{\min}$ follows from the fact that if $X \simeq_{\text{game}} Y$ and $Z \subseteq X$ is stable, then there is a stable $W \subseteq Y$ with $Z \simeq_{\text{game}} W$. Taking $Z$ as all of the bottom chips of $X$, $$\#\text{ of bottom chips of } X = \ell_{\max} (Z) = \ell_{\max} (W) \leq \# \text{ of bottom chips of } Y.$$ Repeating the argument for $Y$ gives the reverse inequality. The rest of the equalities are obvious.)
%\end{proof}

Note that the options $Z \in \operatorname{Opt}(X)$ for which
$\ell_{\min}(Z) = \ell_{\min}(X)-1$ are exactly the options obtained
from removing visible bottom chips, which are the component-end
vertices of $\Gamma(X)$.

Within each connected component, consecutive vertices are adjacent:
otherwise, Lemma \ref{umbrella} would forbid every edge crossing
their gap. Thus, an interior bottom chip has both immediate upper
chips present, while a component-end bottom chip lacks at least one.
Hence, the visible bottom chips are exactly the component ends,
with a singleton contributing one such move.

For a component end $i$, every surviving chip containing $i$ has
$i$ as an endpoint, since all bottom indices of that chip belong
to the same component. The cascade therefore deletes precisely
vertex $i$ and its incident edges, so the resulting option graph
is $\Gamma(X)-i$. Different bottom moves give distinct successor
positions, since each removes its own bottom chip and retains
every other bottom chip.

The associated graph options, counted with multiplicity, form the
multiset $\partial\Gamma(X)$. Equivalently,
\[
\partial\Gamma(X)
=
\left\{\!\left\{
[\Gamma(Z)] :
Z \in \operatorname{Opt}(X),\
\ell_{\min}(Z)=\ell_{\min}(X)-1
\right\}\!\right\}.
\]

\begin{prop}[The Converse]\label{converse}
    For $X, Y$ stable positions,
    $$X \simeq_{\text{game}} Y \implies \Gamma(X) \simeq \Gamma(Y).$$
\end{prop}
\begin{proof}
    Induct on the common chip number $\ell_{\max} (X)$.
    
    The empty case (no chips) is trivially true.

    Suppose $X \simeq_{\text{game}} Y$ and the claim holds for all smaller-chip game-isomorphic positions. Then the option bijection between $\text{Opt}(X)$ and $\text{Opt}(Y)$ restricts to a bijection between the options where a bottom chip was removed. However, by the inductive hypothesis, the graphs of the isomorphic states there are isomorphic. Hence, $\partial \Gamma(X) = \partial \Gamma (Y).$ However, as explained before, we also know that $$|V(\Gamma(X))| = |V(\Gamma(Y))|, \hspace{0.5cm} |E(\Gamma(X))| = |E(\Gamma(Y))|.$$

    Hence,
    $$(|V(\Gamma(X))|, |E(\Gamma(X))|, \partial \Gamma (X)) = (|V(\Gamma(Y))|, |E(\Gamma(Y))|, \partial \Gamma (Y)).$$

    By Theorem \ref{reconstruction}, $\Gamma(X) \simeq \Gamma(Y).$
    
\end{proof}

\section{Formulaic Results}
We now discuss our findings, with some of the enumerative legwork handled in reference.

\subsection{Formula for $K_n$}
This graph-theoretic interpretation of $K_n$ is dependent on the study of the isomorphism equivalence discussed up until now in the paper. We will use the equivalence again in Section 4.2 concerning a derived formula for $K(x)$. 
\begin{thm}\label{enumeration}
    $\forall n\geq0$,\[
K_n=
\sum_{\substack{k_1,\ldots,k_n\ge0\\
\sum_{m=1}^n mk_m\le n}\\}
\prod_{m=1}^n
\binom{c_m+k_m-1}{k_m},
\]

where $c_m$ is the number of connected unit interval graph types on $m$ vertices \cite{SaitohEtAl2010}.
\end{thm}
\begin{proof}

Here, $k_m$ is the number of components with $m$ vertices. Hence, a choice
$(k_1,\ldots,k_n)$ describes a graph with
\[
k_1+2k_2+\cdots+nk_n
\]
vertices. The condition $\displaystyle\sum_{m=1}^n mk_m\le n$ allows us to count every graph of
order at most $n$.

We must choose $k_m$ components from the $c_m$ connected graph types, allowing for repeated
copies and ignoring their order. If the types are
$T_1,\ldots,T_{c_m}$ and $a_i$ is the number of copies of $T_i$, then
\[
a_1+\cdots+a_{c_m}=k_m,\qquad a_i\ge0.
\]
Clearly, the number of choices is
\[
\binom{c_m+k_m-1}{k_m}.
\]

For example, $c_3=2$: the connected 3-vertex types are the path
$P_3$ and the triangle $K_3$. Choosing two such components gives
\[
P_3\sqcup P_3,\qquad
P_3\sqcup K_3,\qquad
K_3\sqcup K_3,
\]
where $\sqcup$ denotes disjoint union. Accordingly, the factor for
$m=3$ and $k_3=2$ is
\[
\binom{2+2-1}{2}=3.
\]

The product combines the independent choices for each component size.
The sum then runs over all possible component counts with total order
at most $n$. Since every graph has a unique decomposition into
connected components, each graph isomorphism class is counted exactly
once. In particular, the choice \[k_1=\cdots=k_n=0\] contributes $1$: the
empty graph.
\end{proof}
\begin{example}
    For $n=3$, the possible component sizes and their contributions are
\[
\begin{array}{c|c}
\text{Component sizes} & \text{\# of choices}\\
\hline
\text{Empty graph} & 1\\
1 & 1\\
1+1 & 1\\
2 & 1\\
1+1+1 & 1\\
2+1 & 1\\
3 & 2
\end{array}
\]
Hence,
\[
K_3=1+1+1+1+1+1+2=8.
\]
\end{example}

\begin{remark}
     In Theorem \ref{enumeration}, replacing the inequality $\displaystyle\sum_{m=1}^n mk_m\le n$ with the corresponding equality $\displaystyle\sum_{m=1}^n mk_m= n$ would instead enumerate $u_n$.
\end{remark}

\subsection{Formula for $K(x)$}
    As promised, we will now apply to proof of isomorphism equivalence to find a rational formula for $K(x)$. Since we have proven both directions of the argument, Theorem \ref{targetthm} is true and we can thus conclude Corollary \ref{form}. Therefore we get the necessary equality in Equation \ref{target}, \[K_n=\sum_{k=0}^nu_k.\]

    Define $K(x)$ as in Equation \ref{gf}, and similarly
    \begin{equation*}
        U(x)=\sum_{n\geq0}u_nx^n,
    \end{equation*}
    with $u_0 = 1$. Thus by convolution of ordinary generating functions, \begin{equation}\label{kinu}
        K(x)=\sum_{n\geq0}K_nx^n=\left(\sum_{j\geq0}u_jx^j\right)\left(\sum_{k\geq0}x^k\right)=\frac{U(x)}{1-x}.
    \end{equation}
    Now we will use known results of types of unit interval graphs on $n$ vertices, most of which were proven by Hanlon in \cite{Hanlon1982}. 
    \begin{definition}
        Recalling Lemma \ref{reversal}, let us define the number of connected unit interval graphs on $n$ vertices up to isomorphism to be $c_n$. This has a known closed form \cite{SaitohEtAl2010}, \[c_n=\frac12\left(C_{n-1}+\binom{n-1}{\lfloor(n-1)/2\rfloor}\right).\]
        Moreover, let \(\Psi(x)\) be $c_n$'s ordinary generating function,
        \[\Psi(x)=\sum_{n\geq1}c_nx^n.\] 
    \end{definition}
    So,\[\Psi(x)=\sum_{n\geq1}\frac12\left(C_{n-1}+\binom{n-1}{\lfloor(n-1)/2\rfloor}\right)x^n=\frac14\left(\frac{1+2x}{\sqrt{1-4x^2}}-\sqrt{1-4x}\right).\]
    First, we have the result from Hanlon (using the Euler transform for similar reasons as in Theorem \ref{enumeration}), \begin{equation}\label{hanlon}
        U(x)=\exp\left(\sum_{k\geq1}\frac{\Psi(x^k)}{k}\right).
    \end{equation}
Therefore, by Equations \ref{kinu} and \ref{hanlon}, \begin{equation}\label{genfun}
        K(x)=\frac1{1-x}\exp\left(\sum_{k\geq1}\frac1{4k}\left(\frac{1+2x^k}{\sqrt{1-4x^{2k}}}-\sqrt{1-4x^k}\right)\right).
    \end{equation}
    Which can be written as the Euler product\[K(x)=\frac{1}{1-x}\prod_{n\geq1}(1-x^n)^{-c_n},\]and provides a rational generating function for $K_n$. Thus we can write it in expanded form as a formal power series:\[K(x)=1+2x+4x^2+8x^3+17x^4+38x^5+93x^6+244x^7+691x^8+2080x^9+\cdots.\] 

    \begin{prop}
        $K(x)$ is not algebraic.
    \end{prop}
    \begin{proof}
        By Equation \ref{genfun}, the analytic continuation has a square-root branch point at \[x=\frac1{\sqrt[n]{4}},\qquad \text{for every odd }n\geq1.\]At each of these points, only the $n$th summand adds that particular singularity, so it is still there after exponentiation. Infinitely many branch points imply $K(x)$ is transcendental \cite{FlajoletSedgewick2009}. 
    \end{proof} 

\section{Conclusion}

\subsection{Summary}
Because the isomorphisms of unit interval graphs are completely characterized, we arrive at a complete characterization of game isomorphisms.

\begin{cor}
    For $X, Y \in \mathcal T_n$, the following are equivalent:
    \begin{enumerate}[label=(\roman*)]
        \item $X \simeq _{\text{game}} Y$
        \item $\Gamma(X) \simeq \Gamma(Y)$
        \item $Y$ is a combination of permutations of the order of supported components of $X$, translations of supported components of $X$, and reversals (flips) of supported components of $X$.
        \item $X, Y$ are isomorphic as posets.
    \end{enumerate}
\end{cor}
\begin{proof}
    (i) $\iff$ (ii) is Theorem \ref{targetthm} and follows from Proposition \ref{forward} and Proposition \ref{converse}.

    (ii) $\iff$ (iii) follows from Lemma \ref{reversal}.

    (iv) $\implies$ (i), using the poset isomorphism on the chips to map each sub-state of $X$ to a sub-state of $Y$. By induction, this is indeed a game isomorphism. Meanwhile, (iii) $\implies$ (iv), since permutations of the order of supported components of $X$, translations of supported components of $X$, and reversals (flips) of supported components of $X$ are poset isomorphisms on the chips.
\end{proof}
\begin{remark}
    Unfortunately, game isomorphism is too simple to help studying the optimal strategy of Pyranim, which we plan to discuss in another paper. However, game isomorphism \textit{does} give a relatively natural sense of isomorphism for Dyck paths (and many more Catalan-type objects).
\end{remark}

We also have the following results:
$$|\mathcal T_n| = C_{n+1},\qquad |\mathcal S_n^*|=K_n-K_{n-1}=u_n,$$
$$K(x) = \frac{U(x)}{1-x},$$
$$K_n=
\sum_{\substack{k_1,\ldots,k_n\ge0\\
\sum_{m=1}^n mk_m\le n}}
\prod_{m=1}^n
\binom{c_m+k_m-1}{k_m}.$$

\begin{remark}
    It is much easier to compute $K_n$ by the above formula than it is to find every isomorphism class of the $C_{n+1}$ game states by Definition \ref{recursive}, which is how $K_n$ was originally computed, up to $n = 13$. You can find more values of $K_n$ in \cite{OEISA400758}.
\end{remark}

\subsection{Other Findings}
During our study, we also found other combinatorial objects that align with a representation of the unique game play of Pyranim---a few of which we have cataloged below. \begin{enumerate}
    \item Type-$A$ positive-root ideals, where the number of remaining chips is equal to the size of the root-ideal.
    \item Nonnesting partitions, where each move splits one maximal interval into several.
    \item 132-avoiding permutations, where each move corresponds to an increase in the number of inversions.
    \item Nondecreasing parking functions, where the game raises parking preferences as chips disappear.
    \item The positive roots of $\mathfrak{sl}_{n+1}$, where chip removal corresponds to adjoining a root vector and closing under Borel action.
\end{enumerate}

\subsection{Open Questions}

The correspondence between Pyranim states and unlabeled unit
interval graphs suggests several natural extensions of the game.
In particular, one may consider variations of the underlying
geometry, support structure, and associated partial order.
We conclude with several questions for further investigation.

\begin{enumerate}
    \item How does the number of game-isomorphism classes change
    when Pyranim is played on a three-dimensional tetrahedron? More generally, can one characterize the number
    of isomorphism classes for a $d$-dimensional pyramid?

    \item Is there a higher-dimensional analogue of the
    correspondence between Pyranim states and unlabeled
    unit interval graphs? In particular, can the isomorphism
    classes of three-dimensional Pyranim be described using
    a natural family of graphs, hypergraphs, or posets?

    \item Can similar games be constructed on other families
    of posets, such as the positive root posets of types
    $B_n$, $C_n$, and $D_n$? If so, do their isomorphism
    classes admit analogous combinatorial interpretations?

    \item What is the expected number of isomorphism classes when considering random starting points of the games---similar to competitive, three-move ballot checkers? Moreover, what is the probability that two randomly selected game states from $\mathcal T_n$ are isomorphic? The study of \[\lim_{n\to\infty}\sqrt[n]{\mathbb P(X\simeq Y \mid X,Y\in\mathcal T_n)}\,\,\stackrel{?}{=}\,1\]is of particular interest, as well as the subexponential decay of this probability.
\end{enumerate}

There are plenty more interesting combinatorial objects to study that rely on the bijective properties of our game to the Catalan numbers and UUIGs. One of these is the directed acyclic graphs seen directly in the formulation of Pyramid Nim in \cite{pyranim}, and other correspondences between these objects are discussed in \cite{pyramid, Kato2024}. Clearly, this game contains some incredible algebraic properties, some of which have yet to be explored.  

\newpage
\printbibliography
 
\end{document}